\documentclass[11pt]{amsart}
\usepackage{amsmath}
\usepackage{amsfonts}

\usepackage{amscd}
\usepackage{amsthm}
\usepackage{amssymb} \usepackage{latexsym}
\usepackage{eufrak}
\usepackage{euscript}
\usepackage{epsfig}
\usepackage{graphics}
\usepackage{array}
\usepackage{enumerate}
\usepackage{dsfont}
\usepackage{color}
\usepackage{wasysym}
\usepackage{hyperref}
\usepackage{pdfsync}

\newcommand{\bel}[1]{\begin{equation}\label{#1}}

\newcommand{\be}{\begin{equation}}

\newcommand{\ba}{\begin{eqnarray}}
\newcommand{\ea}{\end{eqnarray}}

\newcommand{\qe}{\end{equation}}

\newcommand{\Hmm}[1]{\leavevmode{\marginpar{\tiny%
$\hbox to 0mm{\hspace*{-0.5mm}$\leftarrow$\hss}%
\vcenter{\vrule depth 0.1mm height 0.1mm width \the\marginparwidth}%
\hbox to
0mm{\hss$\rightarrow$\hspace*{-0.5mm}}$\\\relax\raggedright #1}}}

\theoremstyle{theorem}
\newtheorem{thm}{Theorem}[section]

\theoremstyle{example}

\theoremstyle{corollary}

\theoremstyle{lemma}
\newtheorem{lem}[thm]{Lemma}
\theoremstyle{definition}
\newtheorem{defi}[thm]{Definition}
\theoremstyle{proof}
\theoremstyle{remark}
\newtheorem{rem}[thm]{Remark}

\begin{document}

\title{One Example about the Relationship between the CD Inequality and CDE' Inequality}
\author{Yijin Gao}
\email{yjgao@ruc.eud.cn}
\address{Department of Mathematics,Information School,
Renmin University of China,
Beijing 100872, China
}

\begin{abstract}In this paper,we will give an easy example to satisfy that we can not conclude CDE' Inequality just from the CD Inequality.
\end{abstract}
\maketitle

\section{Introduction}
The curvature-dimension inequality(CD-inequlaity ) was firstly introduced by Bakry-Emery as a substitute of the lower Ricci curvature bound of the underlying space.It was first studies on finite graphs by Lin-Yau{\cite{9}}{\cite{11}}.

Then some special cases $CD(0,\infty)$ are being studid by Liu-Peyerimhoff.They prove an optimal eigenvalue ratio estimate for finite weighted graphs satisfying the $CD$ inequality.Also,they show taking Cartesian products to be an efficient way for construcing new weighted graphs satisfying $CD(0,\infty)$.{\cite{6}}{\cite{7}}.

In 2015,Paul-Lin-Liu-Yau prove Li-Yau type estimates for bounded and positiv solutions of the heat equation on graphs,under the assumption of the curvature-dimension inequality $CDE'(0,n)$,which can be consider as a notion of curvature for graphs.So the relation between $CD$ and $CDE'$ is becoming more and more important for the study of the ricci estimate on finite graphs.{\cite{2}}

In 2014,2015,Munch introduced a new version of a curvature-dimension inequality for non-negative curvature.He used this inequlaity to prove a logarithmic Li-Yau inequality on finite graphs.The new calculus and the new curvature-dimension inequality coinciden with the common ones.In the case of graphs,they coincide in a limit.In that case,the new curvature-dimension inequality gives a more general concept of curvature on graphs.Then he showed the connection between the $CDE'$ and the $CD \psi$ inequality.Also,he showed that the $CDE'$ inequality implies the $CD$ inequality.{\cite{3}}{\cite{4}}

So can $CD$ inequality implies $CDE'$?It is a problem.In this paper,the author gives a special case to show that the hypothesis is not right in some special cases.

\indent The paper is organized into three parts:\\
\indent Chapter $1$ is the introduction of the graph, the Laplacians and CD,CDE' inequalities on it.\\
\indent Chapter $2$ introduces some basic conclusions and lemmas  in order to get the main result. These conclusions include some definitions such as local finite graph, weighted graph and the calculation of the $\Gamma$ operator.\\
\indent Chapter $3$ is the main conclusion of this thesis.\\

\section{GRAPHS, LAPLACIANS CD INEQUALITIES AND CDE' INEQUALITIES}

\indent Given a graph $G=(V,E)$, for an $x\in V$, if there exists another $y\in V$ that satisfies $(x,y)\in E$, we call them are neighbors, and written as $x\sim y$. If there exists an $x\in V$ satisfying $(x,x)\in E$, we call it a self-loop. For a graph $G=(V,E)$.The neighborhood and the degree of a vertx $x \in V$ are defined,respectively,as $N_x=\{y \in V:xy \in E\}$ and $d_x=|N_x|$.For notational simplicity we work with $\mu_x=1/d_x$.\\
\indent Now we will introduce some basic definitions and theorems before we get the main results.
\begin{defi}{\rm(difference in valuations)}For a function $f \in R^V$ and two vertices$x,y \in V$ denote by$f(x,y)=f(y)-f(x)$ the difference in valuations.\end{defi}
\begin{defi} {\rm (locally finite graph)}\ \ We call a graph G is a locally finite graph if for any $x\in V$, it satisfies $\#\{y\in V|y\sim x\}<\infty$. Moreover, it is called connected if there exists a sequence $\{x_i\}_{i=0}^n$ satisfying: $x=x_0\sim x_1\sim \cdots \sim x_n=y$.\end{defi}
\begin{defi} {\rm (Laplacians on locally finite graphs)}\ \ On a locally finite graph $G=(V,E,\mu)$ the Laplacian has a form as follows: $$\triangle f(x)=\mu_x\sum_{y\in N_x}(f(y)-f(x)),\quad \forall f\in C_0(V).$$\end{defi}
\begin{defi} {\rm (gradient operator $\Gamma$)}\ \ The operator $\Gamma$ is defined as follows:$$\Gamma (f,g)(x)=\frac{1}{2}(\triangle(fg)-f\triangle g-g\triangle f)(x)$$.\end{defi}
\indent Always we write $\Gamma (f,f)$ as $\Gamma(f)$.
\begin{defi} {\rm (gradient operator $\Gamma_{i}$)}\ \ The operator $\Gamma_{i}$ is defined as follows:$$\Gamma_{0}(f,g)=fg$$\\ $$\Gamma_{i+1}(f,g)=\frac{1}{2}(\triangle(\Gamma_{i}(f,g))-\Gamma_{i}(f,\triangle g)-\Gamma_{i}(\triangle f,g))$$.\end{defi}
\indent Also we have $\Gamma_2(f)=\Gamma_2(f,f)=\frac{1}{2}\triangle\Gamma(f)-\Gamma(f,\triangle f)$.

\begin{defi}{\rm ($CD(K,n)$ condition)} \ \ We call a graph satisfies $CD(K,n)$ condition if for any $x \in V$,we have$$\Gamma_2(f)(x)\ge \frac{1}{n}(\triangle f)^2(x)+K \Gamma(f)(x).\quad K\in \mathbb{R}.$$\end{defi}
\begin{defi} {\rm ($CDE(K,n)$ condition)}\ \ Let $f:V\rightarrow\mathbb{R}^+$ satisfy $f(x)>0$, $\triangle f(x)<0$. We call a graph satisfies $CDE(x,K,n)$ condition if for any $x\in V$, we have$$\Gamma_2(f)(x)-\Gamma\left(f,\frac{\Gamma(f)}{f}\right)(x)\geqslant\frac{1}{n}(\triangle f)(x)^2+K\Gamma(f)(x).\quad K\in \mathbb{R}.$$\end{defi}
\begin{defi} {\rm ($CDE'(K,n)$ condition)}\ \ We say that a graph $G$ satisfies the exponential curvature dimension inequality $CDE(K,n)$ if for any positive function $f:V \to R^+$ such that $\Delta f(x)<0$,we have
$$\Gamma_2(f)(x)-\Gamma(f,\frac{\Gamma(f)}{f}(x))\ge \frac{1}{n}f(x)^2(\Delta log f)(x)^2+k\Gamma(f)(x)$$\end{defi}.

\begin{lem} $$\Gamma(f)(x)=\frac{1}{2}\mu_x\sum_{y \in N_x} f(x,y)^2.$$ \\
Here we define $N_x=\{y \in V:xy \in E\}$ and $d_x=|N_x|.$For notational simplicity we work with $\mu_x=\frac{1}{d_x}.$
\end{lem}
{\bf Proof}
\indent We have
\begin{equation*}
\left.
\begin{aligned}
&\Gamma (f)(x)=\frac{1}{2}\triangle(f^2)(x)-f(x)(\triangle f)(x)\\
&=\frac{1}{2}\mu_x \sum_{y \in N_x}(f^2)(x,y)-f(x)\mu_x \sum_{y \in N_x} f(x,y)\\
&=\frac{1}{2}\mu_x\sum_{y \in N_x}(f(x,y)(f(y)+f(x))-2f(x,y)f(x))\\
&=\frac{1}{2}\mu_x\sum_{y \in N_x}f(x,y)^2.\\
\end{aligned}
\right.
\end{equation*}

\begin{lem} $$\Gamma_2(f)(x)=\frac{1}{2}(\triangle f)^2(x)+\mu_x \sum_{y \in N_x}\mu_y\sum_{z \in N_y} (f(y,z)^2-\frac{1}{2}f(x,z)^2)).$$
{\bf Proof}
\indent We have
\begin{equation*}
\left.
\begin{aligned}
&\triangle (\Gamma(f))(x)=\mu_x\sum_{y \in N_x} \Gamma(f)(x,y)
&=\mu_x \sum_{y \in N_x}\frac{1}{2}\mu_y\sum_{z \in N_y}(f(y,z)^2-f(x,y)^2)\\
\end{aligned}
\right.
\end{equation*}
\indent and
\begin{equation*}
\left.
\begin{aligned}
&\Gamma(f,\triangle f)(x)=\frac{1}{2}(\triangle(f\cdot \triangle f)(x)-f(x)\cdot(\triangle^2 f)(x)-(\triangle f)^2(x))\\
&=-\frac{1}{2}(\triangle f)^2(x)+\frac{1}{2}\mu_x \sum_{y \in N_x}((f \triangle f)(x,y)-f(x)(\triangle f)(x,y))\\
&=-\frac{1}{2}(\triangle f)^2(x)+\frac{1}{2}\mu_x \sum_{y \in N_x}f(x,y)(\triangle f)(y)\\
&=-\frac{1}{2}(\triangle f)^2(x)+\frac{1}{2}\mu_x\sum_{y \in N_x}f(x,y)\mu_y \sum_{z \in N_y}f(y,z)\\
\end{aligned}
\right.
\end{equation*}

\indent thus
\begin{equation*}
\left.
\begin{aligned}
&\Gamma_2(f)(x)=\frac{1}{2}\triangle (\Gamma(f))(x)-\Gamma(f,\triangle f)(x)\\
&=\frac{1}{2}(\triangle f)^2(x)+\frac{1}{2}\mu_x\sum_{y \in N_x}\mu_y \sum_{z \in N_y}(\frac{1}{2}f(y,z)^2-\frac{1}{2}f(x,y)^2-f(x,y)f(y,z))\\
&=\frac{1}{2}(\triangle f)^2+\frac{1}{2}\mu_y\sum_{y \in N_x} \mu_y \sum_{z \in N_y}(f(y,z)^2-\frac{1}{2}f(x,z)^2)\\
\end{aligned}
\right.
\end{equation*}

\end{lem}

\begin{lem}If $\Delta f(x)<0$ in $x \in V$,$CDE'(K,n)$ implies $CDE(K,n).$
{\bf Proof}
\indent Let $f:V \to R+$ be a positive function for which $\Delta f(x) <0$.Since $log s \le s-1$for all positive $s$,we can write
$$\Delta log f(x)=\sum_{y \sim x}(log f(y)-log f(x))\le \sum_{y \sim x}\frac{f(y)-f(x)}{f(x)}=\frac{\Delta f(x)}{f(x)} <0.$$
\indent Hence squaring everything reverses the above inequality and we get
$$(\Delta f(x))^2 \le f(x)^2(\Delta log f(x))^2,$$
\indent and thus $CDE(K,n)$ is satisfied
$$\Gamma_2(f)(x)\ge \frac{1}{n}f(x)^2(\Delta log f)(x)^2+k\Gamma(f)(x) > \frac{1}{n}(\Delta f)(x)^2+k\Gamma(f)(x).$$
\end{lem}

\begin{lem}
\indent The $CDE'$ inequality implies $CD$ inequlaity.\\
It was the work in {\cite{3}}
\end{lem}

\section{BASIC CONCLUSION}
\begin{rem} In this section,we just concern the easiest graph: $x$ is the initial point and its neighborhood is $y$,also $y$ has another neighborhood $z$.so the graph consists three points which created from $x$,$x$ and $z$ are not connected.We give the special graph a name:$EG$.Also,we concern the special case in $CD$ and $CDE'$ inequlaity:$CD(0,n)$ and $CDE'(0,n)$.\\
\end{rem}

\begin{lem}
\indent If the $EG$ satisfies $CD(0,m)$,then we have $m \ge 2$.\\
{\bf Proof}
\indent For the simplicity,we rewrite $f(x)=x$,$f(y)=y$ and $f(z)=z$.Then according to the lemma before ,we can have the following:
$\Delta f(x)=y-x$,$\Gamma(f)(x)=\frac{1}{2}(y-x)^2$.
\indent Also,we need to calculate the $\Gamma_2(f)$,
\begin{equation*}
\left.
\begin{aligned}
&\Gamma_2(f)(x)=\frac{1}{2}(y-x)^2+\frac{1}{4}[(z-y)^2-\frac{1}{2}(z-x)^2+(x-y)^2]\\
&=\frac{1}{2}(y-x)^2+\frac{1}{4}(z,y)^2-\frac{1}{8}(z-x)^2+\frac{1}{4}(x-y)^2\\
&=\frac{1}{8}[z^2+z(2x-4y)+8y^2+5x^2-12xy]\\
\end{aligned}
\right.
\end{equation*}
\indent Absolutely,we will use the knowledge of quadratic function.The minimun of $z^2+z(2x-4y)$ can be obtained when $z=2y-x$.Then we take $z=2y-x$ in the equality.
\begin{equation*}
\left.
\begin{aligned}
&8y^2+5x^2-12xy+(2y-x)^2+(2y-x)(2x-4y)\\
&=8y^2+5x^2-12xy+4y^2+x^2-4xy+4xy-8y^2-2x^2+4xy\\
&=4x^2+4y^2-8xy\\
\end{aligned}
\right.
\end{equation*}
\indent So we have the estimate of $\Gamma_2(f)$\\
$$\Gamma_2(f) \ge \frac{1}{2}(x^2+y^2-2xy)=\frac{1}{2}(x-y)^2 \ge \frac{1}{m}(y-x)^2$$
\indent So we have $m \ge 2$.

\end{lem}

\begin{thm}
\indent If the $EG$ satisfies $CD(0,m)$,we can take the special value in the point $y$,then it may not be satisfied in $CDE'(0,m)$.\\
{\bf Proof}
\indent Without loss of generality,we assume $x=1$.According the theorem before,we need to find some special $y$ that make sure the $EG$ does not satisfy $CDE'(0,2)$.\\
\indent As before,the value of $\Gamma_2(f)$ is given
$$\Gamma_2(f)=\frac{3}{4}(x-y)^2+\frac{1}{4}(z,y)^2-\frac{1}{8}(z,x)^2$$
\indent So the next work is to calculate the value of $\Gamma(f,\frac{\Gamma(f)}{f})$
\begin{equation*}
\left.
\begin{aligned}
&\Gamma(f,\frac{\Gamma(f)}{f}=\frac{1}{2}\Delta(\Gamma(f))-\frac{1}{2}f(x)\Delta(\frac{\Gamma(f)}{f(x)})-\frac{1}{2}\Delta(f)\frac{\Gamma(f)}{f(x)}\\
&=\frac{1}{2}\Delta(\Gamma(f))-\frac{1}{2}\Delta(\frac{\Gamma(f)}{f})-\frac{1}{2}\Delta f \Gamma(f)\\
&=I_1-I_2-I_3.\\
\end{aligned}
\right.
\end{equation*}
\indent Then we get the value of $I_1$.
\begin{equation*}
\left.
\begin{aligned}
&I_1=\frac{1}{2}(\Gamma(f)(y)-\Gamma(f)(x))\\
&=\frac{1}{2}{\frac{1}{4}[(y-x)^2+(z,y)^2]-\frac{1}{2}(y-x)^2}\\
&=\frac{1}{2}[\frac{1}{4}(y-x)^2+\frac{1}{4}(y-z)^2-\frac{1}{2}(y-x)^2]\\
&=\frac{1}{8}(y-z)^2-\frac{1}{8}(y-x)^2.\\
\end{aligned}
\right.
\end{equation*}
\indent Also the value of $I_2$.
\begin{equation*}
\left.
\begin{aligned}
&I_2=\frac{1}{2}\frac{\Gamma(f)(y)}{y}-\frac{1}{2}\frac{\Gamma(f)(x)}{x}\\
&=\frac{1}{2}\frac{1/4[(y-x)^2+(z-y)^2]}{y}-\frac{1}{4}(y-x)^2\\
\end{aligned}
\right.
\end{equation*}
\indent Then the value of $I_3$.
\begin{equation*}
\left.
\begin{aligned}
&I_3=\frac{1}{2} \Delta(f) \Gamma(f)(x)\\
&=\frac{1}{2}(y-x)\frac{1}{2}(y-x)^2\\
&=\frac{1}{4}(y-x)^3\\
\end{aligned}
\right.
\end{equation*}

\indent We get the following inequality according the definition of $CDE'(m,0)$
$$\frac{3}{4}(y-x)^2+\frac{1}{4}(y-z)^2-\frac{1}{8}(z-x)^2-\frac{1}{8}(y-z)^2+\frac{1}{8}(y-x)^2+\frac{1}{8y}(y-x)^2+\frac{1}{8y}(y-z)^2\\
-\frac{1}{4}(y-x)^2+\frac{1}{4}(y-x)^3 \ge \frac{1}{2}(log y )^2$$

\indent Firstly,we deal with the polynomial with $z$ like the situation before.
\begin{equation*}
\left.
\begin{aligned}
&\frac{1}{8}(y-z)^2-\frac{1}{8}(z-x)^2+\frac{1}{8y}(y-z)^2\\
&=\frac{1}{8}(y^2+z^2-2yz-z^2-1+2z)+\frac{1}{8y}(y^2+z^2-2yz)\\
&=\frac{1}{8}y^2-\frac{1}{4}yz+\frac{1}{4}z-\frac{1}{8}+\frac{1}{8}y+\frac{1}{8y}z^2-\frac{1}{4}z\\
&=\frac{1}{8y}z^2-\frac{1}{4}yz+\frac{1}{8}y^2+\frac{1}{8}y-\frac{1}{8}\\
\end{aligned}
\right.
\end{equation*}

\indent Here,we take $z=\frac{\frac{1}{4}y}{\frac{1}{4}\frac{1}{y}}=y^2$.Then we get
$$\frac{1}{8y}y^4-\frac{1}{4}y^3+\frac{1}{8}y^2+\frac{1}{8}y-\frac{1}{8}=-\frac{1}{8}y^3+\frac{1}{8}y^2+\frac{1}{8}y-\frac{1}{8}.$$

\indent Then the inequality becomes
\begin{equation*}
\left.
\begin{aligned}
&The ..left..of ..the ..inequality\\
&=\frac{1}{4}(y-1)^3+\frac{5}{8}(y-1)^2+\frac{1}{8y}(y-1)^2-\frac{1}{8}y^3+\frac{1}{8}y^2+\frac{1}{8}y-\frac{1}{8}\\
&=\frac{1}{4}(y^3-3y^2+3y-1)+\frac{5}{8}y^2+\frac{5}{8}-\frac{10}{8}y-\frac{1}{8}y^3+\frac{1}{8}y^2+\frac{1}{4}y-\frac{3}{8}+\frac{1}{8y}\\
&=\frac{1}{8}y^3+y(\frac{3}{4}-\frac{10}{8}+\frac{1}{8}+\frac{1}{8})+\frac{1}{8y} \\
&=\frac{1}{8}y^3-\frac{1}{4}y+\frac{1}{8y}\\
\end{aligned}
\right.
\end{equation*}
\indent So the inequality becomes:
$$y^3-2y+\frac{1}{y} \ge 4(log y)^2$$.
\indent Here we distinguish $y$ into two parts:$y>1$ and $0 <y <1$\\
\indent I:First,we concern the situation $y>1$,because $e^x \ge x+1$,so we have $log y < y-1$\\
\indent Then we need  to prove
$$y^3-2y+\frac{1}{y} \ge 4(y-1)^2$$
\indent It equals
$$y^3-4y^2+6y-4+\frac{1}{y} \ge 0$$
\indent we set the function $h(y)=y^3-4y^2+6y-4+\frac{1}{y}.$\\
\indent The first derivate of $h(y)$ is
$$h'(y)=3y^2-8y+6-\frac{1}{y^2}$$
\indent The second derivate of $h(y)$ is
$$h''(y)=6y-8+2\frac{1}{y^3}$$
\indent Then we use the cauchy inequality to get that
$$6y-8+\frac{1}{y^3}=2y+2y+2y+2\frac{1}{y^3}-8=2(y+y+y+\frac{1}{y^3})-8 \ge 2*4 (yyy \frac{1}{y^3})^(1/4)-8=0$$
\indent So the function $h'(y)$ is increasing when $y >1$.
$$h'(y)>h'(1)=3-8+6-1=0$$
\indent So the function $h(y)$ is increasing too.
$$h(y) \ge h(1)=1-4+6-4+1=0$$.
\indent When $y>1$,we can get the $CDE'(0,m)$
\\

\indent II:Second,we concern the situation$y <1$ as our discussion.we need to proof the following
$$y^3-2y+\frac{1}{y} \ge 4(log f)^2$$
\indent Actually, it equals that
$$y^3-2y+\frac{1}{y} \ge 4(log\frac{1}{y})^2$$
\indent Just as discussed before,then we have to prove
$$y^3-2y+\frac{1}{y} \ge 4(\frac{1}{y}-1)^2$$
\indent Here we set $Q(y)=y^3-2y+\frac{1}{y}-4(\frac{1}{y}-1)^2=y^3-2y+\frac{9}{y}-\frac{4}{y^2}-4$. The next step is to analysis $Q(y)$

\begin{equation*}
\left.
\begin{aligned}
&y^3-2y+\frac{9}{y}-\frac{4}{y^2}-4\\
&=y^3-y^2+y^2-y-y+1-5+\frac{5}{y}+\frac{4}{y}-\frac{4}{y^2}\\
&=(y-1)(y^2+y-1-\frac{5}{y}+\frac{4}{y^2})\\
&=(y-1)(y^2-y+2y-2+1-\frac{1}{y}-\frac{4}{y}+\frac{4}{y^2}\\
&=(y-1)^2(y+2+\frac{1}{y}-\frac{4}{y^2})\\
&=(y-1)^2(y-1+3-\frac{3}{y}+\frac{4}{y}-\frac{4}{y^2})\\
&=(y-1)^3(1+\frac{3}{y}+\frac{4}{y^2})\\
\end{aligned}
\right.
\end{equation*}

\indent Of course $1+\frac{3}{y}+\frac{4}{y^2} >0$,but $y-1 <0$,so $(y-1)^3(1+\frac{3}{y}+\frac{4}{y^2}) <0$,which is contrary to what we need.Clearly,we can set $y=0.1$,then the left is less than the right.So we get the conclusion that we can not conclude $CDE'$ just from $CD$ situation.

\end{thm}

\end{document}